\documentclass[11pt,a4paper]{article}
\usepackage{amsfonts,amscd,amsthm,amsgen,amsmath,amssymb}
\usepackage[all]{xy}
\newtheorem{prop}{Proposition}[section]
\newtheorem{thm}{Theorem}[section]

\newtheorem{lemp}{Lemma}[section]

\theoremstyle{remark}

\theoremstyle{definition}
\newtheorem{defn}{Definition}[section]

\title{Cyclic Higgs bundles and the affine Toda equations}
\author{David Baraglia\footnote{Email: david.baraglia@anu.edu.au}\\
Mathematical sciences institute \\ 
The Australian National University \\
Canberra ACT 0200, Australia}
\date{\today}
\begin{document}
\maketitle
\begin{abstract}
We introduce a class of Higgs bundles called cyclic which lie in the Hitchin component of representations of a compact Riemann surface into the split real form of a simple Lie group. We then prove that such a Higgs bundle is equivalent to a certain class of solutions to the affine Toda equations. We further explain this relationships in terms of lifts of harmonic maps.

\end{abstract}


\section{Introduction}
The Higgs bundle equations and affine Toda field equations lie at the intersection of some intensively studied topics in integrable systems. Namely they are examples of harmonic maps from a surface to a homogeneous space and they are also realized as certain dimensional reductions of the Yang-Mills self-duality equations. It is not so surprising then to find that there is a relationship between the two equations. This relationship has already been observed by Aldrovandi and Falqui \cite{ald} in which it is proved that a certain real form of the affine Toda equations yield solutions to the Higgs bundle equations of a specific form, which in this paper are called cyclic. In this article we will provide a proof of the converse, that any cyclic Higgs bundle arises from a solution to the affine Toda equations. It is worth noting that we are able to derive the Toda equations from the Higgs bundle equations by imposing only a few natural assumptions.\\

There are at least two points of view in which one can see the relationship between Higgs and Toda. On the one hand such equations give rise to flat connections and representations of the fundamental group of a surface. On the other hand they are both related to harmonic maps into certain homogeneous spaces. For the most part we will fix our attention on the flat connection side of the story, leaving some observations about the corresponding harmonic maps to the final section of the paper.\\

There is a vast literature on the topic of harmonic maps from a surface to a homogeneous space. The affine Toda equations and the Higgs bundle equations are two particular cases of this and our objective here is to prove a relation between the two. The Toda equations are related to harmonic maps into a quotient $G/T$ of a simple Lie group $G$ modulo a maximal abelian $T$. Usually $G$ is taken to be the compact real form and $T$ a maximal torus, however for our investigations it is necessary to consider the case where $G$ is the split real form of a simple Lie group and $T$ is generated by a certain Cartan subalgebra. Now while it is true that all maximal tori in a compact simple Lie group are conjugate, for other real forms there generally exists several inequivalent classes of maximal abelian subgroups modulo conjugation. In terms of the Toda equations this translates to the existence of several inequivalent real forms of the equations. The particular real form of the equations of interest here is governed by the connection to Higgs bundles.\\

The Higgs bundle equations are related to harmonic maps into a quotient $G/K$ where $K$ is the maximal compact subgroup of $G$. In attempting to relate the Toda equations to Higgs bundles we are roughly speaking asking when this harmonic map lifts to a harmonic map $G/(K \cap T)$ where $T$ is a maximal abelian subgroup. In the case where $G$ is the split real form of a simple Lie group Hitchin \cite{hit1} identified a component of the moduli space of Higgs bundles which we call the {\em Hitchin component}. For Higgs bundles in the Hitchin component we are able to give a simple characterization of the Higgs bundles arising from the affine Toda equations. A Higgs bundle in the Hitchin component corresponds to a solution of the Toda equations if and only if it is cyclic, where a cyclic Higgs bundle is one in which the Higgs field is cyclic as an element of the Lie algebra in the sense of Kostant \cite{kost}, except at finitely many points where it is allowed to vanish. Alternatively the cyclic Higgs bundles are characterized by the property that they correspond to fixed points of a certain finite cyclic subgroup of the $\mathbb{C}^*$-action on the moduli space of Higgs bundles.\\

In Section \ref{algsec} we begin with a review of some structure theory of simple Lie algebras, the principal three dimensional subgroups and cyclic elements. Following this Section \ref{higsec} presents a quick review of Higgs bundles leading to the construction of Higgs bundles for the Hitchin component. Then in Section \ref{speccase} we then introduce cyclic Higgs bundles and prove that they reduce to certain real forms of the affine Toda equations.

We present a review of the affine Toda equations in Section \ref{sectoda} and their description in terms of flat connections. We then show that under the appropriate reality conditions the affine Toda equations provide a solution to the Higgs bundle equations. A small modification of the Toda equations puts them in a form applicable to any Riemann surface and we then show that the appropriate real form of these equations is equivalent to a cyclic Higgs bundle. Finally in Section \ref{remsec} we take a look at the Higgs Toda relationship on the level of harmonic maps.


\section{Algebraic preliminaries}\label{algsec}
We review some results in the theory of complex simple Lie algebras. We will also take this opportunity to fix some notation that will be used throughout.


\subsection{Principal three-dimensional subalgebras}\label{ptd}
Let $\mathfrak{g}$ be a complex simple Lie algebra of rank $l$. We let $G$ be the adjoint form of $\mathfrak{g}$, that is $G = {\rm Aut}(\mathfrak{g})_0$ is the identity component of the automorphism group of $\mathfrak{g}$. It is well known that the Lie algebra of $G$ is $\mathfrak{g}$ and that the center of $G$ is trivial.

We choose a Cartan subalgebra $\mathfrak{h}$, let $\Delta \subset \mathfrak{h}^*$ be the roots, choose a system of positive roots $\Delta^+$ and corresponding set of simple roots $\Pi$. Fix a corresponding basis $\{h_\beta, e_\alpha,e_{-\alpha} | \; \beta \in \Pi, \; \alpha \in \Delta^+\}$. Here we use the convention that if $\alpha$ is a root then $h_\alpha \in \mathfrak{h}$ is the corresponding coroot defined by
\begin{equation*}
\beta(h_\alpha) = 2\frac{(\alpha , \beta)}{(\alpha , \alpha)}.
\end{equation*}
Kostant \cite{kost} defines a subalgebra (unique up to conjugacy) called the {\em principal $3$-dimensional subalgebra} of $\mathfrak{g}$. We construct the principal $3$-dimensional subalgebra as follows. Let
\begin{equation}\label{x}
x = \dfrac{1}{2}\sum_{\alpha \in \Delta^+}h_\alpha
\end{equation}
then $x = \sum_{\alpha \in \Pi}r_\alpha h_\alpha$ for some positive half-integers $r_\alpha$. We use these to further define
\begin{equation*}
e = \sum_{\alpha \in \Pi}\sqrt{r_\alpha}e_\alpha, \; \; \;
\tilde{e} = \sum_{\alpha \in \Pi}\sqrt{r_\alpha}e_{-\alpha}.
\end{equation*}
Then we define $\mathfrak{s}$ as the linear span of $\{x,e,\tilde{e}\}$. We must verify that $\mathfrak{s}$ is a subalgebra:
\begin{lemp}\cite{Oni1}
For any $\beta \in \Pi$ we have $\beta(x)=1$.
\begin{proof}
Let $R_\beta$ denote the reflection in the hyperplane in $\mathfrak{h}^*$ orthogonal to $\beta$, namely
\begin{equation*}
R_\beta(\alpha) = \alpha - 2\frac{(\alpha,\beta)}{(\beta,\beta)}\beta.
\end{equation*}
The dual action of $R_\beta$ on $\mathfrak{h}$ is then
\begin{equation}\label{reflection}
R_\beta^t (h) = h-\beta(h)h_\beta.
\end{equation}
From this one verifies the relation $R_\beta^t(h_\alpha) = h_{R_\beta(\alpha)}$. We then have that
\begin{equation*}
R_\beta^t(x) = \dfrac{1}{2}\sum_{\alpha \in \Delta^+} h_{R_\beta(\alpha)}.
\end{equation*}
However we also have that for $\alpha \in \Delta^+$, $R_\beta(\alpha) \in \Delta^+$ except for $\alpha = \beta$ in which case $R_\beta(\beta) = -\beta$. It follows that $R_\beta^t(x) = x - h_\beta$, hence by (\ref{reflection}), $\beta(x)=1$.
\end{proof}
\end{lemp}

Given a root $\lambda \in \Delta$, let $\lambda = \sum_{\alpha \in \Pi}n_\alpha \alpha$. Then we define the {\em $\Pi$-height} of $\lambda$ to be the integer ${\rm height}(\lambda) = \sum_{\alpha \in \Pi}n_{\alpha}$. The above lemma shows that for $y \in \mathfrak{g}_\alpha$, we have $[x,y] = {\rm height}(\alpha)y$, that is $x$ is the grading element corresponding to the gradation of $\mathfrak{g}$ by height. We now deduce the following commutation relations for $\mathfrak{s}$:
\begin{equation*}
[x,e] = e, \; \; \; [x,\tilde{e}] = -\tilde{e}, \; \; \; [e,\tilde{e}] = x.
\end{equation*}
Thus $\mathfrak{s}$ is a copy of $\mathfrak{sl}(2,\mathbb{C})$.\\

We have that the element $e$ is regular, that is it has an $l$-dimensional centralizer spanned by elements $e_1, \dots , e_l$. Moreover on restriction to $\mathfrak{s}$ the adjoint representation decomposes into irreducible subspaces:
\begin{equation}\label{gdec}
\mathfrak{g} = \bigoplus_{i=1}^l V_i.
\end{equation}
We can take $e_1, \dots , e_l$ as highest weight elements of $V_1, \dots , V_l$ which shows there are indeed $l$ summands. Since $\mathfrak{s}$ itself must appear as one of the $V_i$, we take it to be $V_1$ so we may take $e_1 = e$. Let $m_1, \dots , m_l$ denote the exponents of $\mathfrak{g}$. These can be described as follows \cite{Oni2}: arrange the positive roots of $\mathfrak{g}$ into an array with the $k$-th row consisting of all roots of height $k$, filling in rows from right to left. Then the lengths of the columns from left to right are the exponents. In particular we see that if the highest root $\delta$ has height $M$ then $m_l = M$ is the largest exponent.

Returning to the decomposition of $\mathfrak{g}$ given by (\ref{gdec}), we have that the dimensions of $V_1, \dots , V_l$ are $(2m_1+1), \dots , (2m_l + 1)$, from which we may write
\begin{equation}\label{gdecomp}
\mathfrak{g} = \bigoplus_{i=1}^l {\rm S}^{2m_i}(V)
\end{equation}
where $V$ is the $2$-dimensional fundamental representation for $\mathfrak{s}$. Observe that $[x,e_l] = Me_l$ and thus $e_l$ is a highest weight vector.

We may also decompose $\mathfrak{g}$ according to the action of $x$, i.e. by height:
\begin{equation*}
\mathfrak{g} = \bigoplus_{m= -M}^M \mathfrak{g}_m.
\end{equation*}


\subsection{Coxeter and cyclic elements}
The key algebraic feature underpinning much of what follows is presence of a certain distinguished element of $G$:
\begin{equation*}
g = {\rm exp}(2\pi i x / h)
\end{equation*}
where $h = M+1$ is the Coxeter number. We have that $g^h = 1$ and the eigenvalues of ${\rm Ad}_g$ are the $h$-th roots of unity. The highest weight vectors $e_1, \dots , e_l$ for the representation of the principal three dimensional subgroup on $\mathfrak{g}$ span a Cartan subalgebra. Kostant \cite{kost} shows that $g$ acts on this Cartan subalgebra as a Coxeter element for the Weyl group. As such we will call $g$ a {\em Coxeter element} of $G$. Conversely every lift to $G$ of a Coxeter element of a Cartan subalgebra is of this form.

We can decompose $\mathfrak{g}$ into the eigenspaces of ${\rm Ad}_g$:
\begin{equation*}
\mathfrak{g} = \bigoplus_{m=0}^M \mathfrak{g}^g_m,
\end{equation*}
where $\mathfrak{g}^{g}_m$ is the $e^{2\pi i m/ h}$-eigenspace of ${\rm Ad}_g$. In particular $\mathfrak{g}^g_0 = \mathfrak{h}$ is the Cartan subalgebra and $\mathfrak{g}^g_1$ is the sum of the root spaces of the simple roots together with the root space of the lowest root $-\delta$, that is
\begin{equation*}
\mathfrak{g}^g_1 = \bigoplus_{\alpha \in \Pi \cup \{-\delta\}} \mathfrak{g}_\alpha.
\end{equation*}

It is clear that $g$ treats the simple roots and the lowest root on an equal footing. We let $\overline{\Pi} = \Pi \cup \{-\delta\}$.\\

Let $X \in \mathfrak{g}_1^g$. We say that $X$ is {\em cyclic} if $X = \sum_{\alpha \in \overline{\Pi}} c_\alpha e_\alpha$ where the constants $c_\alpha$ are non-zero for all $\alpha \in \overline{\Pi}$. More generally, we say that $Y \in \mathfrak{g}$ is cyclic if it is $Y$ conjugate to such an $X$ under the adjoint action of $G$. Cyclic elements have a simple description in terms of the invariant polynomials on $\mathfrak{g}$.

Recall that an invariant polynomial on $\mathfrak{g}$ is an element of the polynomial algebra $S^*(\mathfrak{g}^*) = \bigoplus_{i=0}^\infty S^k(\mathfrak{g}^*)$ that is invariant under the natural action of $G$. There are homogeneous generators $p_1 , \dots , p_l$ for the invariant polynomials on $\mathfrak{g}$ such that for all elements $f \in \mathfrak{g}$ of the form
\begin{equation*}
f = \tilde{e} + \alpha_1 e_1 + \dots + \alpha_l e_l,
\end{equation*}
we have
\begin{equation}\label{poly}
p_i(f) = \alpha_i.
\end{equation}
We also have that $p_i$ has degree $m_i+1$. It turns out that an element $X \in \mathfrak{g}$ is cyclic if and only if $p_1(X) = p_2(X) = \dots = p_{l-1}(X) = 0$, and $p_l(X) \neq 0$ \cite{kost}. Put another way $X$ is cyclic if every invariant polynomial of degree less than $h$ without constant term vanishes on $X$ but there is an invariant polynomial without constant term that is non-vanishing on $X$.\\

Let $X,Y \in \mathfrak{g}^g_1$ be cyclic. There exists an element $t \in T$ such that ${\rm Ad}_t X = \lambda Y$ for some non zero constant $\lambda$. Moreover $X$ and $Y$ are conjugate under $T$ if and only if $p_l(X) = p_l(Y)$.


\section{Higgs bundles and the Hitchin component}\label{higsec}


\subsection{Higgs bundles}
Higgs bundles were introduced by Hitchin in \cite{hit2} after considering dimensional reduction of the Yang-Mills self-duality equations to a surface. A {\em Higgs bundle} over a Riemann surface $\Sigma$ is a pair $(E,\Phi)$ where $E$ is a holomorphic vector bundle on $\Sigma$ and $\Phi$ is a holomorphic section of ${\rm End}(E) \otimes K$, where $K$ is the canonical bundle on $\Sigma$. We call $\Phi$ a {\em Higgs field}. The {\em slope} of a Higgs bundle $(E,\Phi)$ is the slope of the bundle $E$, namely $\mu(E) = {\rm deg}(E)/{\rm rank}(E)$. There is a $\mathbb{C}^*$-action on Higgs bundles which acts by multiplication of the Higgs field by a non-zero complex number.\\

Closely related to Higgs bundles are the {\em Higgs bundle equations} for a pair $(\nabla_A , \Phi)$ on a Hermitian vector bundle $E$ over a Riemann surface. Here $\nabla_A$ is a unitary connection on $E$ and $\Phi$ is a section of ${\rm End}(E) \otimes K$. The Higgs bundle equations are:
\begin{eqnarray}
F_A + \left[ \Phi , \Phi^* \right] &=& 0 \label{hbe1} \\
\nabla_A^{0,1}\Phi &=& 0 \label{hbe2}
\end{eqnarray}
where $F_A$ is the curvature of $\nabla_A$ and $\nabla_A^{0,1}$ is the $(0,1)$ part of $\nabla_A$.\\

The relation between Higgs bundles and the Higgs bundle equations is as follows: if $(\nabla_A , \Phi)$ is a solution of the Higgs bundle equations then $\nabla_A^{0,1}$ defines a holomorphic structure on $E$ for which $\Phi$ is holomorphic and this defines a Higgs bundle pair $(E,\Phi)$. Conversely suppose $(E,\Phi)$ is a Higgs bundle. If we choose a Hermitian form $h$ on $E$ then we get a Chern connection $\nabla_A$, the unique connection such that $\nabla_A(h)=0$ and $\nabla_A^{0,1}$ defines the holomorphic structure on $E$. Then one asks if we can choose the Hermitian form $h$ such that $(\nabla_A , \Phi)$ solves the Higgs bundle equations. Even with the obvious restriction that $E$ must have degree $0$ it is not always possible to find such a metric. The result of Hitchin \cite{hit2} and Simpson \cite{simp1} is that such a metric exists if and only if the Higgs bundle $(E,\Phi)$ is {\em polystable}. In this case the unitary connection $\nabla_A$ is unique. We now explain the polystable condition.\\

If $(E,\Phi)$ is a Higgs bundle then a subbundle $F \subset E$ is called {\em $\Phi$-invariant} if $\Phi(F) \subseteq F \otimes K$. The Higgs bundle $(E,\Phi)$ is said to be {\em stable} if for every proper, non-zero $\Phi$-invariant holomorphic subbundle $F \subset E$ we have the inequality $\mu(F) < \mu(E)$. If this condition is relaxed to the non-strict inequality $\mu(F) \le \mu(E)$ we say $E$ is {\em semistable}. We also say that $(E,\Phi)$ is {\em polystable} if it is a direct sum of stable Higgs bundles of the same slope. \\

The above discussion generalizes to the case of principal $G$-bundles where $G$ is a complex Lie group with Lie algebra $\mathfrak{g}$. A $G$-Higgs bundle is a pair $(P,\Phi)$ where $P$ is a principal $G$-bundle over $\Sigma$ and $\Phi$ is a holomorphic section of $\mathfrak{g}^P \otimes K$, where $\mathfrak{g}^P$ is the adjoint bundle associated to $P$. 

Similarly the Higgs bundle equations have a generalization to principal $G$-bundles. We will assume that $G$ is semisimple. Then for a principal $G$-bundle the equivalent of a Hermitian form is a reduction of structure to the maximal compact subgroup $G^\rho$. The adjoint bundle $\mathfrak{g}^\rho$ is a real subbundle of the complex adjoint bundle $\mathfrak{g}^P$ and defines an antilinear involution $\rho : \mathfrak{g}^P \to \mathfrak{g}^P$. If $x \otimes y \mapsto k(x,y)$ is the Killing form on $\mathfrak{g}^P$ then corresponding to $\rho$ is the Hermitian form $x \otimes y \mapsto -k(x , \rho(y))$. The Hermitian adjoint of ${\rm ad}_x$ is $-{\rm ad}_{\rho(x)}$, hence we will denote $-\rho(x)$ by $x^*$. 

It is now straightforward to generalize the Higgs bundle equations. Suppose $P$ is a principal $G$-bundle with a reduction of structure $P_H \subset P$ to a principal $H$-subbundle, with $H = G^\rho$ the maximal compact subgroup. Then the Higgs bundle equations for a connection $\nabla_A$ on $P_H$ and a section $\Phi$ of $\mathfrak{g}^P \otimes K$ are given by equations (\ref{hbe1}), (\ref{hbe2}).\\

Note also that if $(\nabla_A , \Phi)$ is a solution of the Higgs bundle equations for a principal $G$-bundle then $\nabla = \nabla_A + \Phi + \Phi^*$ is a flat $G$-connection and is determined by its monodromy representation $\pi_1(\Sigma) \to G$. Thus Higgs bundles on $\Sigma$ are closely related to representations of the fundamental group of $\Sigma$.


\subsection{The Hitchin component}\label{hitcomp}

Higgs bundles provide a geometric interpretation for representations of the fundamental group of a surface into a simple Lie group. A particular explicit construction of Higgs bundles by Hitchin \cite{hit1} identifies a component of the space of representations into the split real form of a Lie group. This component is often called the Hitchin component. We consider within the Hitchin component those representations for which only the highest holomorphic differential is non-vanishing and show in this case the Higgs bundle equation reduces to a set of Toda equations.\\

We review the construction of Higgs bundles in \cite{hit1}. Throughout we will assume $\Sigma$ has genus $g > 1$. Recall that if ${\bf M}$ is the moduli space of Higgs bundles $(E,\Phi)$ for the group $G$ then there is a map
\begin{equation*}
p : {\bf M} \to \bigoplus_{i = i}^l H^0(\Sigma , K^{m_i+1})
\end{equation*}
obtained by applying the invariant polynomials $p_1, \dots , p_l$ to $\Phi$. We will construct a section of this map.\\

Consider the bundle
\begin{equation}\label{E}
E = \bigoplus_{m=-M}^M \mathfrak{g}_m \otimes K^m
\end{equation}
where $\mathfrak{g}_m$ is the $m$-eigenspace of the grading element $x$. This is the adjoint bundle of $\mathfrak{g}$ associated to the principal $G$-bundle $P = P_1 \times_i G$ where $P_1$ is the holomorphic principal ${\rm SL}(2,\mathbb{C})$-bundle associated to $K^{1/2}\oplus K^{-1/2}$ (for a choice of $K^{1/2}$) and $i: {\rm SL}(2,\mathbb{C}) \to G$ is the inclusion corresponding to the principal $3$-dimensional subalgebra $\mathfrak{s}$. We see from (\ref{gdecomp}) that $E$ is independent of choice of $K^{1/2}$ and defines for us a holomorphic bundle with the structure of $\mathfrak{g}$.

Next we construct a Higgs field $\Phi$. Let $q_1, \dots , q_l$ be holomorphic differentials of degrees $m_1+1, \dots , m_l +1$. Then we may define $\Phi \in H^0(\Sigma , E \otimes K)$ as follows:
\begin{equation*}
\Phi = \tilde{e} + q_1 e_1 + \dots + q_l e_l.
\end{equation*}
Here, since $\tilde{e} \in \mathfrak{g}_{-1}$, we regard $\tilde{e}$ as a section of $(\mathfrak{g}_{-1} \otimes K^{-1}) \otimes K$ and similarly $q_ie_i$ can be considered as a section of $(\mathfrak{g}_{i} \otimes K^i ) \otimes K$ so that $\Phi$ is a well defined holomorphic section of $E \otimes K$. The Higgs bundles $(E,\Phi)$ constructed here are polystable and in fact correspond to smooth points of the moduli space ${\bf M}$ \cite{hit1}. We have from (\ref{poly}) that $p_i(\Phi) = q_i$ and it follows that $(q_1 , \dots , q_l ) \mapsto (E,\Phi)$ is our desired section.\\

So far we have constructed Higgs bundles with holonomy in $G$. We show that in fact the Higgs bundles constructed have holonomy in the split real form of $G$. We define an involutive automorphism $\sigma$ on $\mathfrak{g}$ as follows:
\begin{equation}\label{sig}
\sigma(e_i) = -e_i, \; \; \; \sigma(\tilde{e}) = -\tilde{e}.
\end{equation}
Since $\mathfrak{g}$ is obtained from $e_1, \dots , e_l$ by repeated application of ${\rm ad}_{\tilde{e}}$, we see that such a $\sigma$ must be unique. We state some properties of $\sigma$ \cite{hit1}; $\sigma$ exists and defines a Cartan involution corresponding to the split real form of $\mathfrak{g}$. That is if $\rho$ is a compact anti-involution defining a Hermitian metric solving the Higgs bundle equations then in fact $\rho$ and $\sigma$ commute and $\lambda = \rho \sigma$ is an anti-involution for the split real form of $\mathfrak{g}$. Clearly our Higgs field satisfies $\sigma \Phi = -\Phi$. If $(\nabla_A , \Phi)$ satisfies the Higgs bundle equations then so does $(\sigma^*\nabla_A , \sigma^*\Phi)$ and hence so does $(\sigma^*\nabla_A , \Phi)$. Uniqueness of $\nabla_A$ now shows that $\sigma$ is preserved by $\nabla_A$. Note also that $\lambda \Phi = - \rho (\Phi) = \Phi^*$. It now follows that the flat connection $\nabla_A + \Phi + \Phi^*$ has holonomy in the split real form of $G$. Therefore we have constructed representations of the fundamental group of $\Sigma$ into the split real form.\\

Our section $s : \bigoplus_{i = i}^l H^0(\Sigma , K^{m_i+1}) \to {\bf M}$ takes values in the smooth points of ${\bf M}$ and identifies the vector space $V = \bigoplus_{i = i}^l H^0(\Sigma , K^{m_i+1})$ with a submanifold of ${\bf M}$. Moreover the map $s$ actually maps $V$ into the moduli space ${\bf M}^\lambda$ of flat $G^\lambda$-connections where $G^\lambda$ is the split real form of $G$. Around the smooth points ${\bf M}^\lambda$ is a manifold of dimension $(2g-2) \dim(G^\lambda)$. Moreover we find from the Riemann-Roch theorem that $\dim(V) = (2g-2)\dim(G^\lambda)$ as well. Since the image $s(V) \subset {\bf M}^\lambda$ is open, closed and connected it defines a smooth component of the space of representations of the fundamental group in $G^\lambda$. We call this the {\em Hitchin component}.


\section{Cyclic Higgs bundles}\label{speccase}
We now consider a special case of this construction in which all but the highest differential are set to zero. We will call such a Higgs bundle {\em cyclic}. Thus our Higgs field is of the form

\begin{equation}\label{higgs}
\Phi = \tilde{e} + qe_l.
\end{equation}
Here $q \in {\rm H}^0(\Sigma,K^{M+1})$ is a holomorphic $(M+1)$-differential, where $M = m_l$. Thus our Higgs field is cyclic as an element of $\mathfrak{g}$ at all points where $q \ne 0$.\\

Solutions of the form (\ref{higgs}) have the property that they are fixed points in the moduli space of Higgs bundles of the action of certain roots of unity in $\mathbb{C}^*$. From the properties of the element $x$ we have
\begin{equation*}
[x,\Phi] = -\tilde{e} + Mqe_l.
\end{equation*}
Let $g = {\rm exp}(2\pi i x/(M+1))$ be the Coxeter element. It follows that
\begin{equation*}
{\rm Ad}_g \Phi = {\rm exp}(2\pi i M/(M+1))\Phi = \omega\Phi
\end{equation*}
where $\omega$ is an $(M+1)$-th root of unity. Now if $(\nabla_A,\Phi)$ is the solution to the Higgs bundle equation corresponding to $\Phi$ then we may gauge transform this solution by $g$ to get the equivalent solution $(g^*\nabla_A, \omega\Phi)$, we can then use the $U(1)$-action to obtain the solution $(g^*\nabla_A,\Phi)$. The uniqueness theorem for solutions to the Higgs bundle equations now implies that $A$ is gauge invariant under $g$, which in turn is equivalent to covariant constancy of $g$ with respect to the connection $\nabla_A$. Thus the action of $g$ and $\nabla_A$ on sections of $E$ commutes. Observe that the Cartan subalgebra $\mathfrak{h}$ is the unity eigenspace of ${\rm Ad}_g$ from which it follows that $\nabla_A$ preserves the subbundle of $E$ corresponding to $\mathfrak{h}$ and thus the connection form $A$ is $\mathfrak{h}$-valued.\\

We now seek an explicit form of the Higgs bundle equations for the special class of solutions (\ref{higgs}). For this we adopt the following point of view: given the holomorphic data $(E,\Phi)$, we seek a Hermitian metric $H$ on $E$ such that the Chern connection $A = H^{-1}\partial H$ satisfies the Higgs bundle equations:
\begin{equation*}
F_A + [\Phi,\Phi^*] = 0
\end{equation*}
where $F_A = \overline{\partial}(H^{-1}\partial H)$ is the curvature of $\nabla_A$.

More specifically, in our case a Hermitian metric is a reduction of the structure group of $E$ to the maximal compact subgroup $G^\rho$ of $G$. This is equivalent to finding an anti-linear involution $\rho : E \to E$ preserving the Lie algebra structure of the fibres, i.e. it is a choice of compact real form of $\mathfrak{g}$ on each fibre of $E$. The associated Hermitian form is then
\begin{equation*}
H(u,v) = h_\rho(u,v) = -k(u,\rho(v))
\end{equation*}
where $k$ is the Killing form on $\mathfrak{g}$. Recall that for $u \in \mathfrak{g}$, the adjoint of $u$ thought of as the endomorphism ${\rm ad}_u$ is given by $x^* = -\rho(x)$. Now since $G$ is the adjoint form of $\mathfrak{g}$ we can think of $G$ as a group of linear transformations in $\mathfrak{g}$. Corresponding to the real form $\rho$ on $\mathfrak{g}$ we may define a map $\rho : G \to G$ by $\rho(g) = (g^*)^{-1}$. Here the adjoint is with respect to the corresponding Hermitian form $H$. Note that the differential of $\rho$ at the identity is the anti-involution $\rho$ which justifies denoting both by $\rho$.\\

In order to write down an arbitrary such anti-involution $\rho$ we first define the following fixed anti-involution $\hat{\rho}$ of $\mathfrak{g}$:
\begin{equation}\label{compact}
\hat{\rho}(h_\alpha) = -h_\alpha, \; \; \; \hat{\rho}(e_\alpha) = -e_{-\alpha}, \; \; \; \hat{\rho}(e_{-\alpha}) = -e_\alpha.
\end{equation}
The fixed set of $\hat{\rho}$ defines a compact real form of $\mathfrak{g}$ and moreover any other compact real form  is conjugate to $\hat{\rho}$ by an inner automorphism.

We emphasize that $\rho$ is an anti-involution on the adjoint bundle $E$ while $\hat{\rho}$ is a fixed anti-involution on the Lie algebra $\mathfrak{g}$. Consider now a trivialisation of $E$ over an open subset $U$ which identifies the fibres of $E$ with $\mathfrak{g}$ preserving the Lie algebra structure. The anti-involution $\rho : E \to E$ defines on each fibre $E_x$ an anti-involution $\rho_x : E_x \to E_x$ and upon identifying $E_x$ with $\mathfrak{g}$, $\rho_x$ defines an anti-involution on $\mathfrak{g}$. Hence for each $x \in U$ there exists a $p(x) \in G$ such that $\rho_x = \alpha_x \hat{\rho} \alpha_x^{-1}$ where $\alpha_x = {\rm Ad}_{p(x)}$. Allowing $x$ to vary over $U$ we have a map $p : U \to G$ which can be chosen smoothly such that in the trivialisation over $U$, $\rho = {\rm Ad}_p \circ \hat{\rho} \circ {\rm Ad}_{p}^{-1}$. Thus
\begin{equation}\label{relate}
\rho = {\rm Ad}_h \circ \hat{\rho}
\end{equation}
where $h = p\hat{\rho}(p)^{-1}$. Equation (\ref{relate}) is a local formula relating $\rho$ to $\hat{\rho}$. Note also that ${\rm Ad}_h$ is a positive symmetric operator with respect to $\hat{\rho}$. The following lemma will show that we then have $h = {\rm exp}(2\Omega)$ for some $\Omega \in \mathfrak{g}$ with the property $\hat{\rho}(\Omega) = -\Omega$.\\

\begin{lemp}
Let $\mathfrak{g}$ be a complex Lie algebra with Hermitian form. If $h \in {\rm Aut}(\mathfrak{g})$ is positive and symmetric then $h = e^H$ where $H$ is symmetric and $H \in {\rm der}(\mathfrak{g})$ is a derivation of $\mathfrak{g}$.
\begin{proof}
We may decompose $\mathfrak{g}$ into the eigenspaces of $h$
\begin{equation*}
\mathfrak{g} = \bigoplus_\lambda \mathfrak{g}_\lambda.
\end{equation*}
Since $h$ is an automorphism we have $\left[ \mathfrak{g}_{\lambda_1} , \mathfrak{g}_{\lambda_2} \right] \subseteq \mathfrak{g}_{\lambda_1 \lambda_2}$. Now we may write $h = e^H$ where $H = {\rm log}(h)$ is symmetric. Further $H$ acts on $\mathfrak{g}_\lambda$ with eigenvalue ${\rm log}(\lambda)$. The fact that $H$ is a derivation follows immediately since ${\rm log}(\lambda_1\lambda_2) = {\rm log}(\lambda_1) + {\rm log}(\lambda_2)$.
\end{proof}
\end{lemp}

Consider how to find the Chern connection associated to the Hermitian metric $h_{\rho}$. In a local trivialisation where $\rho = {\rm Ad}_h \circ \hat{\rho}$ we have:
\begin{eqnarray*}
h_{\rho}(u,v) &=& -k(u,{\rm Ad}_{h}(\hat{\rho}(v))) \\
&=&  h_{\hat{\rho}}(u,{\rm Ad}_{h^{-1}} v).
\end{eqnarray*}
We can regard $h_{\hat{\rho}}$ as a fixed Hermitian form on $\mathfrak{g}$, hence we can think of the metric $h_\rho$ as being locally associated to the matrix $h^{-1} = {\rm exp}(-2{\rm ad}_\Omega)$. The connection form in this trivialisation is then $A = h \partial h^{-1}$.\\

We will need the following:
\begin{lemp}\label{uniq}
The anti-involution $\rho$ solving the Higgs bundle equations is unique
\begin{proof}
Any other solution is related by a $G$-valued gauge transformation $L: E \to E$ such that $L$ preserves the holomorphic structure and $\Phi$. We will argue that $L$ is the identity. It is a general fact \cite{hit2} that a gauge transformation that is holomorphic and commutes with $\Phi$ is in fact covariantly constant with respect to the unitary connection $\nabla_A$. Now as in (\ref{E}) we may write
\begin{equation*}
E = \bigoplus_{m=-M}^M \mathfrak{g}_m \otimes K^m.
\end{equation*}
But we have seen that $\nabla_A$ is $\mathfrak{h} = \mathfrak{g}_0$-valued, hence $\nabla_A$ preserves the subbundles $E_m = \mathfrak{g}_m \otimes K^m$. If we decompose $L$ into endomorphisms $L_{ij} : E_i \to E_j$ then each $L_{ij}$ must be covariantly constant. But since $E_i \otimes E_j^* = \mathfrak{g}_i \otimes \mathfrak{g}_j^* \otimes K^{i-j}$ there can only be a non-vanishing covariantly constant section if $i=j$ (since $\Sigma$ is assumed to have genus $g>1$). Therefore $L$ preserves $\mathfrak{g}_i$ for all $i$. Since $L$ preserves the Cartan subalgebra $\mathfrak{h} = \mathfrak{g}_0$ then $L$ is valued in the normaliser $N(T)$ of the maximal torus $T$. But $L$ also preserves $\mathfrak{g}_1 = \bigoplus_{\alpha \in \Pi} \mathfrak{g}_{\alpha}$. Therefore the induced action of $L$ on the roots of $\mathfrak{g}$ preserves the choice of simple roots, hence $L$ is valued in the maximal torus $T$.

We also have that $L$ commutes with $\Phi = \tilde{e} + qe_l$, in particular ${\rm Ad}_L \tilde{e} = \tilde{e}$. It follows that $L$ is the identity.
\end{proof}
\end{lemp}

We can now prove the following:
\begin{prop}\label{rhorhohat}
The anti-involution $\rho$ locally has the form $\rho = {\rm Ad}_h \circ \hat{\rho}$ where $h = e^{2\Omega}$ and $\Omega$ is valued in $\mathfrak{h}$.
\begin{proof}
We have seen that the connection $A$ is valued in $\mathfrak{h}$ and that $\nabla_A$ preserves the element $g$. Note that $g$ commutes with the action of $x$ so that $g$ is a well defined global gauge transformation. Consider a new globally defined anti-involution given by $\mu = {\rm Ad}_g \circ \rho \circ {\rm Ad}_{g^{-1}} = {\rm Ad}_q \circ \rho$ where $q = g\rho(g)^{-1}$.

Now since both $g$ and $H$ are covariantly constant we see that $\mu$ yields a metric compatible with $A$. Thus the Chern connection corresponding to $\mu$ is also $A$. Now if we take the Higgs bundle equation for $\rho$:
\begin{equation*}
F_A + [\Phi , -\rho(\Phi)] = 0
\end{equation*}
and apply ${\rm Ad}_g$, noting that $g$ commutes with elements of $\mathfrak{h}$ we get
\begin{equation*}
F_A + [\omega\Phi, -{\rm Ad}_g \rho(\Phi)] = 0
\end{equation*}
but since $\omega{\rm Ad}_g \circ \rho (\Phi) = {\rm Ad}_g \circ \rho \circ {\rm Ad}_{g^{-1}} \Phi = \mu(\Phi)$ we get
\begin{equation*}
F_A + [\Phi,-\mu(\Phi)] = 0.
\end{equation*}
By uniqueness of solutions to the Higgs bundle equations we have that $\mu = \rho$, and so ${\rm Ad}_q$ is the identity. Thus $q$ lies in the centre of $G$ which is trivial since $G$ is the adjoint form, hence $\rho(g) = g$.\\

Locally we write $\rho = {\rm Ad}_h \circ \hat{\rho}$ where $h = e^{2\Omega}$. Then $g = \rho(g) = {\rm Ad}_h (\hat{\rho}(g)) = {\rm Ad}_h(g)$. Thus $g$ and $h$ commute. This implies that $h$ preserves the eigenspaces of $g$, in particular $h$ preserves the Cartan subalgebra $\mathfrak{h}$. Thus $h$ is valued in the normaliser $N(T)$. We will further argue that $h$ is valued in the maximal torus $T$.

Since $N(T)/T$ is a finite group (the Weyl group) there is some positive integer $a$ such that $h^a = e^{2a\Omega}$ is valued in $T$, say $h^a = e^{X}$ where $X$ is valued in $\mathfrak{h}$. However $h^a$ is a positive symmetric operator in the sense that $h^a = h^{a/2}\hat{\rho}(h^{-a/2}) = h^{a/2}(h^{a/2})^*$. Therefore $h^{2a} = e^{X - \hat{\rho}(X)}$ is also positive and symmetric. So $h^{2a}$ has a unique positive symmetric $2a$-th root which is $h$. Taking this root yields
\begin{equation*}
h = e^{(X - \hat{\rho}(X))/2a}.
\end{equation*}
Replacing $2\Omega$ by $(X - \hat{\rho}(X))/2a$ if necessary we have that $h = e^{2\Omega}$ is valued in $T$ and $\Omega$ is valued in $\mathfrak{h}$ with $\hat{\rho}(\Omega) = -\Omega$.
\end{proof}
\end{prop}

Having established that $\Omega$ is $\mathfrak{h}$-valued, we have that $\Omega$ and $\partial \Omega$ commute and the connection is simply $A = -2\partial \Omega$. We have in turn that the curvature is $F_A = -2\overline{\partial}\partial \Omega$.

Next we determine $\Phi^* = -\rho(\Phi)$. Recall that $\Phi = \tilde{e} + qe_l$. Thus $\hat{\rho}(\Phi) = -e - \overline{q}\tilde{e_l}$ where $\tilde{e_l}$ is a lowest weight vector for $\mathfrak{g}$. Thus
\begin{eqnarray*}
\Phi^* &=& {\rm Ad}_h (e+ \overline{q}\tilde{e_l}) \\
&=& \sum_{\alpha \in \Pi} \sqrt{r_\alpha}e^{2\alpha(\Omega)}e_\alpha +\overline{q}e^{-2\delta(\Omega)}\tilde{e_l}.
\end{eqnarray*}
We can now determine the commutator $[\Phi,\Phi^*]$:
\begin{eqnarray*}
[\Phi,\Phi^*] &=& \left[ \sum_{\alpha \in \Pi}\sqrt{r_\alpha}e_{-\alpha} + qe_l ,\sum_{\alpha \in \Pi} \sqrt{r_\alpha}e^{2\alpha(\Omega)}e_\alpha +\overline{q}e^{-2\delta(\Omega)}\tilde{e_l} \right] \\
&=& -\sum_{\alpha \in \Pi}r_\alpha e^{2\alpha(\Omega)}h_\alpha -q\overline{q}e^{-2\delta(\Omega)}h_{-\delta}.
\end{eqnarray*}
Although it has been suppressed, this should really be multiplied by a $dz \wedge d\overline{z}$ term. The Higgs bundle equation for $\Omega$ now becomes:
\begin{equation}\label{toda}
-2\Omega_{z\overline{z}} + \sum_{\alpha \in \Pi}r_\alpha e^{2\alpha(\Omega)}h_\alpha +q\overline{q}e^{-2\delta(\Omega)}h_{-\delta} = 0.
\end{equation}

If $\mathfrak{g}$ is a simple Lie algebra of type $L_n$ then these are a version of the {\em affine Toda field equations} \cite{mik} for the Affine Dynkin diagram $L_n^{(1)}$ which we introduce more generally in Section \ref{sectoda}. Note that since $\Omega$ is $\mathfrak{h}$-valued we have $\rho(\Omega) = \hat{\rho}(\Omega) = -\Omega$ and similarly for $h_\alpha, h_{-\delta}$. Thus $\Omega$ is a real linear combination of the $h_\alpha$. However there is an additional constraint on $\Omega$ that we have yet to consider. By examining this condition we will find that in some cases the above equations will reduce to affine Toda field equations for a smaller affine Dynkin diagram.


\subsection{Symmetry conditions}\label{reality}
The additional constraint on $\Omega$ that we have yet to consider relates to a symmetry property of $\Omega$. Recall that we define an involutive automorphism $\sigma$ on $\mathfrak{g}$ by the properties
\begin{equation*}
\sigma(e_i) = -e_i, \; \; \; \sigma(\tilde{e}) = -\tilde{e}.
\end{equation*}
We also note that $\sigma$ commutes with $\hat{\rho}$ and $\hat{\lambda} = \sigma \hat{\rho} = \hat{\rho}\sigma$ is an anti-involution corresponding to the split real form of $\mathfrak{g}$. As before we have $\sigma(\Phi) = -\Phi$ and $\sigma$ is covariantly constant with respect to $\nabla_A$.\\

We shall argue that in fact $\sigma(\Omega) = \Omega $. First note that since $\sigma A = A$ we have $\sigma(\Omega_{z\overline{z}}) = \Omega_{z\overline{z}}$. Thus the remaining terms of the Toda equations (\ref{toda}) must also be fixed by $\sigma$. We also know \cite{hit1} that $\sigma = \phi \circ \hat{\nu}$ where $\phi$ is the inner automorphism corresponding to a rotation by $\pi$ in the principal $3$-dimensional subgroup and $\hat{\nu}$ is a lift of an automorphism $\nu$ of the Dynkin diagram for $\Pi$. In fact $\nu$ is trivial except for the simple Lie algebras of type $A_n,D_{2n+1},E_6$ for which it has order $2$.

In any case, $\phi$ is a conjugation by an element of $\mathfrak{h}$, so that $\sigma$ agrees with $\hat{\nu}$ on $\mathfrak{h}$, that is we have
\begin{equation*}
\sigma({h_\alpha}) = h_{\nu(\alpha)}, \; \; \; \sigma(h_{-\delta}) = h_{-\delta}.
\end{equation*}
We also have that $\sigma(x) = x$ from which it follows that $r_{\nu(\alpha)} = r_{\alpha}$. By invariance of the right hand side of (\ref{toda}) under $\sigma$ we see that $e^{2\alpha(\Omega)} = e^{2\nu(\alpha)(\Omega)}$, now since $\Omega$ is real with respect to the basis $\{h_\alpha\}$ we may take logarithms to conclude that $\alpha(\Omega) = \nu(\alpha)(\Omega)$ for all $\alpha \in \Pi$, that is $\sigma (\Omega) = \Omega$.\\

From the relation $\sigma(\Omega) = \Omega$ we have that $\rho$ and $\sigma$ commute and thus $\lambda = \sigma \rho = \rho \sigma$ is a real structure on $E$ conjugate to $\hat{\lambda}$, hence $\lambda$ defines a reduction of structure to the split real form $G^\lambda$ of $G$.\\

Now we have that $A$ is invariant under $\rho$ and $\sigma$, hence $\lambda$ also. Therefore $\nabla_A$ has holonomy in $K = G^\rho \cap G^\lambda$, the maximal compact subgroup of $G^\lambda$, in fact the holonomy reduces to a maximal torus $T$ of $K$.

Let $\mathfrak{t}$ denote the Lie algebra of $T$, that is $\mathfrak{t} = \mathfrak{h}^\lambda$. The inclusion $i: \mathfrak{t} \to \mathfrak{h}$ induces the restriction map $r: \mathfrak{h}^* \to \mathfrak{t}^*$. Since $\Omega$ is $\mathfrak{t}$-valued we wish to reinterpret the Toda field equations in terms of the restricted roots on $\mathfrak{t}^*$. Intuitively speaking we are dividing out the action of $\nu$.\\

We can realize $\mathfrak{t}^*$ as a subspace of $\mathfrak{h}^*$. In fact the kernel of the restriction map is the $-1$-eigenspace of $\hat{\nu}^t$ so that $\mathfrak{t}^*$ identifies with the $+1$ eigenspace. Under this identification the restriction map becomes the orthogonal projection map
\begin{equation*}
r(\alpha) = \tfrac{1}{2}(\alpha + \hat{\nu}^t(\alpha)).
\end{equation*}
For any $\beta \in \mathfrak{t}^*$ define the dual vector $\tilde{h}_\beta$ by the usual formula
\begin{equation*}
\alpha(\tilde{h}_\beta) = 2\frac{(\alpha,\beta)}{(\beta,\beta)}
\end{equation*}
for all $\alpha \in \mathfrak{t}^*$. Then, for possibly different (non-zero) constants $\tilde{r}_\beta$ the Toda equation (\ref{toda}) reduces to:
\begin{equation}\label{toda2}
-2\Omega_{z\overline{z}} + \sum_{\beta \in r(\Pi)}\tilde{r}_\beta e^{2\beta(\Omega)}\tilde{h}_\beta + q\overline{q}e^{-2\delta(\Omega)}\tilde{h}_{-\delta} = 0.
\end{equation}
The elements of $r(\Pi)$ together with $-\delta = -r(\delta)$ correspond to an affine Dynkin diagram as explained in Section \ref{sectoda}. That is we have a system of vectors in a Euclidean vector space $\mathfrak{t}^*$ whose Dynkin diagram is an affine Dynkin diagram. Equation (\ref{toda2}) is then the affine Toda field equation for the corresponding affine Dynkin diagram.

All that remains is to actually determine which affine root system they correspond to. If $\nu$ is trivial then $\mathfrak{t} = \mathfrak{h}$ and the affine diagram is just the extended Dynkin diagram $L_n^{(1)}$ obtained by adding to the root system $L_n$ the lowest weight $-\delta$. The remaining cases follow from a straightforward calculation and are as follows: for $\mathfrak{g}$ of type $A_{2n},A_{2n-1},D_{2n+1},E_6$, the corresponding affine diagrams are $A_{2n}^{(2)},C_n^{(1)},B_{2n}^{(1)},F_4^{(1)}$ respectively. Here we are using the Kac names for the affine Dynkin diagrams \cite{carter}.


\section{Geometry of the affine Toda equations}\label{sectoda}


\subsection{The affine Toda equations}

The Toda equations have been extensively studied from the point of view of integrable systems \cite{mik}, \cite{man}, \cite{oli}, \cite{lez} as well as their relation to minimal surfaces and harmonic maps \cite{bolw} \cite{dol1}, \cite{bol}, \cite{dol2}, \cite{dol3}. This should not be surprising given the link between Higgs bundles and harmonic maps. However our treatment of the affine Toda equations allows for more general real forms of the equations than are usually considered.\\

Let $A$ be an $n \times n$ indecomposable generalized Cartan matrix of affine type. We label the rows and columns by $0,1, \cdots , l$, where $l = n-1$. We arrange so that the $l \times l$ matrix $A^0$ formed by removing row $0$ and column $0$ is the Cartan matrix for the root system of a complex simple Lie algebra. So there is a real $l$-dimensional Euclidean vector space $(\mathfrak{h}_{\mathbb{R}},\langle \, , \, \rangle)$ with basis $\{h_i\}$, $i = 1, \ldots , l$ such that if we define a corresponding basis $\{\alpha_i \}$, $i = 1, \ldots , l$ of $\mathfrak{h}_{\mathbb{R}}^*$ by $\alpha_i = 2h_i/\langle h_i , h_i \rangle $ then $A_{ij} = \alpha_j(h_i)$ for $i,j = 1, \ldots , l$.

Moreover \cite{carter} there exists positive integers $a_0,a_1, \ldots , a_l$ with no common factor such that $A (a_0 , \ldots , a_l )^t = 0$. Similarly there exists positive integers $c_0, c_1 , \ldots , c_l$ with no common factor such that $(c_0 , \ldots , c_l)A = 0$. Let us define $h_0 \in \mathfrak{h}_{\mathbb{R}}$ and $\alpha_0 \in \mathfrak{h}_{\mathbb{R}}^*$ by
\begin{eqnarray*}
h_0 &=& -\sum_{i=1}^l c_i h_i, \\
\alpha_0 &=& -\sum_{i=1}^l a_i \alpha_i.
\end{eqnarray*}
Then it follows that
\begin{equation*}
\alpha_i (h_j) = A_{ji}, \, i,j = 0,1,\ldots , l.
\end{equation*}

Now we may define the affine Toda equations
\begin{defn}
Given an $n \times n$ indecomposable generalized Cartan matrix of affine type let $\mathfrak{h}_{\mathbb{R}}$, $\{h_i\}$, $\{\alpha_i\}$ be as above. Let $U$ be an open subset of $\mathbb{C}$. A map $\Omega : U \to \mathfrak{h} = \mathfrak{h}_{\mathbb{R}} \otimes \mathbb{C}$ is said to satisfy the {\em affine Toda equations} if it is a solution to
\begin{equation}\label{tod}
2\Omega_{z\overline{z}} = \sum_{i=0}^l k_i e^{2\alpha_i(\Omega)}h_i,
\end{equation}
where $k_i \in \mathbb{C}$ are constants, $k_i \ne 0, \, i = 0, \ldots , l$.
\end{defn}


\subsection{From affine Toda to Higgs bundles}
It is well known that the affine Toda equations can be re-written as a zero curvature condition for a certain associated connection. The fact that such a formulation exists is intimately connected to their integrability. Alternatively the affine Toda equations admit a description in terms of a class of harmonic maps from a surface into a homogeneous space of the form $G/T$ where $T$ is a maximal torus or more generally the group corresponding to a Cartan subalgebra \cite{bol},\cite{dol3}. In turn, both of these descriptions are closely related to the Higgs bundle equations.\\

We will first re-write the Toda equations as the equations for a flat connection. After imposing appropriate reality conditions we will see that the flat connection decomposes into a unitary connection and Higgs field satisfying the Higgs bundle equations. We note that if different reality conditions are imposed one can arrive at similar equations except involving connections that preserve an indefinite Hermitian form.\\

For now we restrict to the case of an untwisted affine Dynkin diagram. In fact the twisted cases all occur as special cases of the untwisted equations \cite{oli} and we have already seen in Subsection \ref{reality} how this works in some cases.\\

Thus we take a complex simple Lie algebra $\mathfrak{g}$ of rank $l$, we let $\mathfrak{h}$ be the real subspace of a Cartan subalgebra spanned by the coroots $h_1 , \dots , h_l$ corresponding to the simple roots $\alpha_1 , \dots , \alpha_l$. The untwisted affine root system corresponding to $\mathfrak{g}$ amounts to letting $\alpha_0 = -\delta$ be the lowest root and $h_0$ the corresponding coroot. Define $x \in \mathfrak{h}$ by (\ref{x}).\\

Now for each root $\alpha_i$, $i=0,1, \dots , l$, let $e_{\alpha_i}$, be a basis element for the spaces of $\alpha_i$. We can choose a basis such that
\begin{equation*}
[ e_{\alpha_i} , e_{-\alpha_j} ] = \delta_{ij} h_i.
\end{equation*}

Consider the connection $\nabla = \nabla^{1,0} + \nabla^{0,1}$ where
\begin{eqnarray}
\nabla^{1,0} &=& (\partial_{z} - \Omega_z)dz + \Phi \label{std1} \\
\nabla^{0,1} &=& (\partial_{\overline{z}} +\Omega_{\overline{z}})d\overline{z} + \Psi \label{std2}
\end{eqnarray}
and $\Phi$, $\Psi$ are defined by
\begin{eqnarray}
\Phi &=& {\rm Ad}_{e^{-\Omega}}\left( \sum_{i=0}^l c_i e_{-\alpha_i} \right) dz \label{std3} \\
\Psi &=& {\rm Ad}_{e^{+\Omega}}\left( \sum_{i=0}^l d_i e_{+\alpha_i} \right) d\overline{z} \label{std4}
\end{eqnarray}
where the $c_i,d_i$ are constants such that $c_i d_i = k_i$. The different choices of the constants $c_i,d_i$ actually define a $(\mathbb{C}^*)^{l+1}$-family of flat connections. We see that the Toda equations (\ref{tod}) are equivalent to the vanishing of the curvature of $\nabla$. We will call this gauge in which $\nabla$ is written the {\em Toda gauge}.\\

One can write down a whole family of gauge equivalent connections. For any smooth map $H : \Sigma \to \mathfrak{h}$ the gauge transformation by ${\rm Ad}_{e^H}$ changes $\nabla$ to $\nabla'^{1,0} + \nabla'^{0,1}$ where
\begin{eqnarray*}
\nabla'^{1,0} = (\partial_{z} - \Omega_z-H_z)dz + \Phi' \\
\nabla'^{0,1} = (\partial_{\overline{z}} +\Omega_{\overline{z}}-H_{\overline{z}})d\overline{z} + \Psi'
\end{eqnarray*}
with $\Phi'$, $\Psi'$ given by
\begin{eqnarray*}
\Phi' &=& {\rm Ad}_{e^{-\Omega+H}}\left( \sum_{i=0}^l c_i e_{-\alpha_i} \right)dz \\
\Psi' &=& {\rm Ad}_{e^{\Omega +H}}\left( \sum_{i=0}^l d_i e_{\alpha_i} \right)d\overline{z}.
\end{eqnarray*}

Choosing $H$ to be constant accounts for a $(\mathbb{C}^*)^l$ subgroup of the aforementioned $(\mathbb{C}^*)^{l+1}$-family of flat connections. Modulo this gauge symmetry There is a remaining $\mathbb{C}^*$-action which amounts to transformations $\Phi \mapsto \zeta \Phi$, $\Psi \mapsto \zeta^{-1} \Psi$ for $\zeta \in \mathbb{C}^*$. The choice of gauge transformation with $H = \Omega$ will be called the {\em Higgs gauge} for reasons that will soon become clear.\\

Having arrived at the desired zero cuvrature form of the equations, the next step is to impose a reality condition. Let $\hat{\rho}$ be the anti-linear involution of $\mathfrak{g}$ defined as in (\ref{compact}). Note that $\hat{\rho}$ defines the compact form of $\mathfrak{g}$ and that $\hat{\rho}$ preserves $\mathfrak{h}$. We impose the reality condition $\hat{\rho}(\Omega) = -\Omega$, or equivalently $\Omega$ is valued in $\mathfrak{h}_{\mathbb{R}}$, the $-1$ eigenspace of $\hat{\rho}$ on $\mathfrak{h}$. Now to consider a particular real form of the Toda equations we also fix the values of the coefficients $k_i$, $i = 0, \dots , l$ to the following convenient values: $k_i = r_i$ for $i = 1, \dots , l$ where the $r_i$ are defined by $x = \sum_{i=1}^l r_i h_i$ and set $k_0 = 1$. Let us further take $c_i = d_i = \sqrt{r_i}$ for $i = 0, \dots , l$.

Under these conditions the connection $\nabla_A = (\partial_{z} - \Omega_z) + (\partial_{\overline{z}} +\Omega_{\overline{z}})$ is a unitary connection on the trivial bundle $\mathfrak{g} \times U \to U$ with respect to the Hermitian form $h_{\hat{\rho}}(x,y) = -k(x, \hat{\rho}(y))$, where $k$ is the Killing form on $\mathfrak{g}$. Moreover we observe that $\hat{\rho}(\Phi) = -\Psi$, or simply $\Psi = \Phi^*$. It follows easily that:
\begin{prop}
The pair $(\nabla_A , \Phi)$ is a solution of the Higgs bundle equations if and only if $\Omega$ is a solution of the affine Toda equations.
\end{prop}
The above proposition is also proved in \cite[Proposition 3.1]{ald}.\\

Note that in the Higgs gauge $\nabla_A^{0,1} = \overline{\partial}$, that is the Higgs gauge is a holomorphic trivialization of the holomorphic structure determined by $\nabla_A$. Let us also note that in the Higgs gauge, the Hermitian form preserved by $\nabla_A$ becomes
\begin{eqnarray*}
H(x,y) &=& -k(e^{-\Omega} x , \hat{\rho}(e^{-\Omega}y)) \\
&=& -k( x , e^{\Omega} \hat{\rho}( e^{-\Omega} y)) \\
&=& -k (x , e^{2\Omega} \hat{\rho}(y)) \\
&=& -k( x , \rho(y))
\end{eqnarray*}
where $\rho = {\rm Ad}_{e^{2\Omega}} \circ \hat{\rho}$. It is no coincidence that this is the same relation as in Proposition \ref{rhorhohat}.


\subsection{Globalization of the affine Toda equations}
We are interested in putting the affine Toda equations into more geometric terms, applicable to any Riemann surface. Since we have already written the affine Toda equations in terms of the Higgs bundle equations, the remaining step is to introduce transition maps to get a connection on a non-trivial bundle. To achieve this we will also need the data of a holomorphic differential of degree $h = M+1$. This will ultimately lead us back to the class of cyclic Higgs bundles.\\

Let $\Sigma$ be a Riemann surface and $\{ U_i \}$ and open cover by coordinate charts with local holomorphic coordinate $z_i$. On the intersection $U_{ij} = U_i \cap U_j$ let $dz_i = g_{ij}dz_j$ where $g_{ij}$ is a non-vanishing holomorphic function on $U_{ij}$. Further define real valued functions $f_{ij}$ by $g_{ij}\overline{g_{ij}} = e^{2f_{ij}}$.\\

Suppose we have a collection of functions $\Omega_i : U_i \to \mathfrak{h}$ which are related on the overlap $U_{ij}$ by
\begin{equation*}
\Omega_i = \Omega_j + f_{ij}x.
\end{equation*}
So the $\Omega_i$ define a section $\Omega$ of a certain affine bundle over $\Sigma$. We call such a section $\Omega$ a {\em Toda field}. The point of this definition is that for each simple root $\alpha$ we have $\alpha(x)=1$ and hence the $e^{2\alpha(\Omega_i)}dz_i \overline{dz_i}$ define a section of $K \overline{K}$, where $K$ is the canonical bundle of $\Sigma$. On the other hand for the lowest root $-\delta$ we have that the $e^{-2\delta(\Omega_i)}dz_i \overline{dz_i}$ define a section of $K^{-M}\overline{K}^{-M}$. To make sense of the affine Toda equation we introduce a holomorphic degree $(M+1)$-differential $q \in H^0(\Sigma , K^{M+1})$. We may then say that the pair $(\Omega,q)$ consisting of a Toda field and holomorphic differential satisfy the affine Toda equations if
\begin{equation}\label{toda2}
2\Omega_{z\overline{z}} = \sum_{i=1}^l k_i e^{2\alpha_i(\Omega)}h_i + k_0 q \overline{q} e^{-2 \delta(\Omega)} h_0.
\end{equation}
Note that both sides of this equation define sections of $K\overline{K} \otimes \mathfrak{h}$ so this is a well-defined equation on $\Sigma$. If $q = 0$ the equation is locally equivalent to the non-affine Toda equations. On the other hand in a neighborhood of a point where $q \neq 0$ we can find a local coordinate such that $q = (dz)^{M+1}$ and the equation reduces to the affine Toda equations as in (\ref{tod}).\\

On each open subset $U_i$ we have a trivial bundle $\mathfrak{g} \times U_i$ and a pair $(\nabla_{A,i},\Phi_i)$ which in the Higgs gauge have the form

\begin{eqnarray*}
\nabla_{A,i}^{1,0} &=& (\partial_{z_i} - 2{(\Omega_i)}_{z_i})dz_i \\
\nabla_{A,i}^{0,1} &=& \partial_{\overline{z_i}} d\overline{z_i}
\end{eqnarray*}
and
\begin{eqnarray*}
\Phi_i &=&  \left( \sum_{i=1}^l c_i e_{-\alpha_i} + q_i e_\delta  \right) dz_i \\
\Phi_i^* &=& {\rm Ad}_{e^{2(\Omega)_i}}\left( \sum_{i=1}^l d_i e_{\alpha_i} + \overline{q_i} e_{-\delta}  \right) d\overline{z_i}
\end{eqnarray*}
where $q = q_i (dz_i)^{M+1}$ in the local coordinate $z_i$.\\

Next we claim that on the overlaps $\,U_{ij}\,$ the pairs $(\nabla_{A,i},\Phi_i)$ and $(\nabla_{A,j},\Phi_j)$ are gauge equivalent by some gauge transformation $\phi_{ij}$. The $\phi_{ij}$ will satisfy a cocyle condition so we can patch together to get a single Higgs bundle $(\nabla_A , \Phi)$ defined on $\Sigma$. 

Let us take transition maps that preserve the root decomposition $\mathfrak{g} = \mathfrak{h} \oplus \bigoplus_{\alpha \in \Delta} \mathfrak{g}_\alpha$. More specifically we take transitions $\phi_{ij}$ that act trivially on $\mathfrak{h}$ and on a given root space $\mathfrak{g}_\alpha$ acts as multiplication by $g_{ij}^{\alpha(x)}$. Thus a section $X_i$ of $\mathfrak{g}_\alpha \times U_i$ is identified with a section $X_j$ of $\mathfrak{g}_\alpha \times U_j$ over $U_{ij}$ by the relation $X_i = g_{ij}^{\alpha(x)} X_j$. 

It follows that the bundles $\mathfrak{g} \times U_i$ patch together to form precisely the bundle $E$ given by (\ref{E}). Moreover, one easily checks that the Higgs bundles $(\nabla_{A,i},\Phi_i)$ patch together under these transitions to form a Higgs bundle $(\nabla_A , \Phi)$ on $E$. It is immediate from the construction that the pair $(\nabla_A , \Phi)$ is a cyclic Higgs bundle. If we now assume $\Sigma$ is compact of genus $g>1$ then as an immediate corollary we find that $\sigma(\Omega) = \Omega$, where $\sigma$ is the Cartan involution as defined in (\ref{sig}). Combined with the results in Section \ref{speccase} we have thus obtained:
\begin{thm}\label{bigthm}
Let $\Sigma$ be a compact Riemann surface of genus $g>1$. Then any pair $(\Omega,q)$ consisting of a Toda field and holomorphic differential $q$ of degree $(M+1)$ satisfying the affine Toda equations (\ref{toda2}) subject to $k_i = r_i$  for $i=1, \dots , l$, $k_0 = 1$ and $\hat{\rho}(\Omega) = -\Omega$ corresponds to a cyclic Higgs bundle. Conversely any cyclic Higgs bundle gives such a solution of the Toda field equations.
\end{thm}
Moreover, we know that the cyclic Higgs bundles on a compact Riemann surface of genus $g>1$ are parametrized by the corresponding holomorphic differential $q \in H^0(\Sigma , K^{M+1})$ so we also have
\begin{thm}
Let $\Sigma$ be a compact Riemann surface of genus $g>1$. For any holomorphic differential $q \in H^0(\Sigma , K^{M+1})$ there is a unique Toda field $\Omega$ such that $\hat{\rho}(\Omega) = -\Omega$ and $(\Omega,q)$ satisfies the affine Toda equations (with coefficients $k_i = r_i$  for $i=1, \dots , l$, $k_0 = 1$.)
\end{thm}
The uniqueness follows from Lemma \ref{uniq}.


\section{Remark on harmonic maps}\label{remsec}
We have mentioned that the Higgs bundle and Toda equations are related not only to flat connections but also to harmonic maps. We demonstrated the relation between Higgs and Toda using flat connection descriptions. Here we remark on the relation in terms of harmonic maps.\\

Suppose $P_G \to \Sigma$ is a principal $G$-bundle with flat connection $\omega$. Let $\tilde{\Sigma}$ be the universal cover of $\Sigma$, which we view as a right principal $\pi_1(\Sigma)$-bundle over $\Sigma$. The pull-back of $P_G$ to a principal bundle over $\tilde{\Sigma}$ admits global sections in which $\omega$ is trivial. Any two such trivializing sections differ by the right action of $G$. A choice of such section $s$ determines a monodromy represention $\theta : \pi_1(\Sigma) \to G$ of $\omega$ by the relation $s(p\gamma) = s(p) \theta(\gamma)$ for all $\gamma \in \pi_1(\Sigma)$. If we change the section $s$ then $\theta$ changes by conjugation.

Now suppose in addition we are given a reduction of structure to a principal $H$-subbundle $P_H \to P_G$ where $H$ is a subgroup of $G$. Note that we are {\em not} assuming the connection $\omega$ reduces to a connection on $P_H$. With respect to a trivializing section $s$ we can find a map $\tau : \tilde{\Sigma} \to G/H$ such that the pull-back of $P_H$ to $\tilde{\Sigma}$ is given by $\{ s(p)g \, | \, p \in \tilde{\Sigma}, \; g \in \tau(p) \, \}$. If we change the section $s$ to $sg$ for some $g \in G$ then we must change $\tau$ to $g^{-1}\tau$. We also observe that the map $\tau : \tilde{\Sigma} \to G/H$ is $\theta$-equivariant meaning $\tau(p\gamma) = \theta(\gamma)^{-1} \tau(p)$ for all $\gamma \in \pi_1(\Sigma)$.\\

Conversely a representation $\theta : \pi_1(\Sigma) \to G$ determines a flat $G$-connection and a $\theta$-equivariant map $\tau : \tilde{\Sigma} \to G/H$ corresponds precisely to a reduction of structure to $H$ of the associated flat principal $G$-bundle.\\

In the case of Higgs bundles we have a flat $G$-connection, where we assume $G$ is simple and a reduction of structure to the maximal compact subgroup $K \subseteq G$. If $\nabla$ is any flat $G$-connection and we are given a reduction of structure $\tau : \tilde{\Sigma} \to G/K$ to the maximal compact then we may uniquely decompose $\nabla = \nabla_A + \Phi + \Phi^*$ where $\nabla_A$ is a unitary connection, $\Phi$ is a $(1,0)$-valued section of the adjoint bundle and $\Phi$ is the adjoint with respect to the Hermitian form determined by the reduction of structure. A key observation is that the pair $(\nabla_A , \Phi)$ satisfies the Higgs bundle equations if and only if the corresponding map $\tau : \tilde{\Sigma} \to G/K$ is harmonic. This is the link between harmonic maps and Higgs bundles.\\

Similarly in the case of the affine Toda equations, we wrote the equations in terms of a flat $G$-connection for a complex Lie group and a reduction of structure to the subgroup $T$ corresponding to a Cartan subalgebra. Actually, in the case we considered there was even further reduction of structure, but we will overlook this. So there is a reduction of structure map $\tau: \tilde{\Sigma} \to G/T$. When we have a solution of the Toda equations, this map is harmonic. There is a partial converse - a special class of harmonic maps $\tilde{\Sigma} \to G/T$ defined in \cite{bol} called $\tau$-primitive correspond to solutions of the affine Toda equations.\\

Now consider the case of cyclic Higgs bundles. Let $G$ be a complex simple Lie group in adjoint form, $G^{\hat{\rho}}$ the maximal compact, $G^{\hat{\lambda}}$ the split real form, and $K = G^{\hat{\lambda}} \cap G^{\hat{\rho}}$ the maximal compact subgroup of $G^{\hat{\lambda}}$. A cyclic Higgs bundle determines a harmonic map $\tau_{{\rm Higgs}} : \tilde{\Sigma} \to G^{\hat{\lambda}}/K$. Now let $T$ be the subgroup generated by a Cartan subalgebra $\mathfrak{h}$. The compact and split real forms have corresponding anti-involutions $\hat{\rho}$ and $\hat{\lambda}$ which according to the choices we use preserve $\mathfrak{h}$. So there are corresponding abelian subgroups $T^{\hat{\rho}} \subset G^{\hat{\rho}}$, $T^{\hat{\lambda}}
 \subset G^{\hat{\lambda}}$ and $T_K = T^{\hat{\rho}} \cap T^{\hat{\lambda}} \subset K$. The Toda equations imply that the reduction of structure map $\tau_{{\rm Toda}} : \tilde{\Sigma} \to G^{\hat{\lambda}}/ T^{\hat{\lambda}}$ is harmonic. Moreover the reductions of structure $\tau_{{\rm Higgs}},\tau_{{\rm Toda}}$ are compatible in the sense that we have a commutative diagram
\begin{equation*}\xymatrix{
\tilde{\Sigma} \ar[r]^{\tau_{{\rm Higgs}}} \ar[d]_{\tau_{{\rm Toda}}} \ar[dr]^{\tau} & G^{\hat{\lambda}}/K  \\
G^{\hat{\lambda}}/ T^{\hat{\lambda}} & G^{\hat{\lambda}}/(K \cap T^{\hat{\lambda}}) \ar[l] \ar[u]
}
\end{equation*}
It is a straightforward computation to verify that the map $\tau$ is also harmonic. Thus in terms of harmonic maps the relationship between the Higgs bundle and affine Toda equations is as simple as one could hope.

\addcontentsline{toc}{chapter}{Bibliography}

\end{document}